\documentclass[reqno,12pt]{amsart}
\usepackage{graphicx}
\usepackage{epsfig}
\usepackage{amssymb,amsthm,amsmath,amsfonts,amstext,mathrsfs}
\usepackage{latexsym}
\usepackage{url}
\usepackage[all]{xy}
\newtheorem{thm}{Theorem}

\input xy
\xyoption{all}

\def\d{\,{\rm{d}}}

\title[Moments of question mark function]
{Asymptotic formula for the moments of Minkowski question mark function in the interval $[0,1]$.}
\author[Giedrius Alkauskas]{Giedrius Alkauskas}
\begin{document}

\begin{abstract} In this paper we prove the asymptotic formula for the moments of
Minkowski question mark function, which describes the distribution of rationals in the Farey tree.
The main idea is to demonstrate that certain a variation of a Laplace method is applicable in this
problem, hence the task reduces to a number of technical calculations.
\end{abstract}
\maketitle
\begin{center}
\end{center}
\begin{center}
\end{center}
\section{Introduction and main result} The Minkowski question mark function $?(x)$ is defined for
real numbers in the interval $[0,1]$ via the expression
\begin{eqnarray*}
?([a_{0},a_{1},a_{2},a_{3},...])=2(1-2^{-a_{0}}+2^{-(a_{0}+a_{1})}-2^{-(a_{0}+a_{1}+a_{2})}+...),
\label{min}
\end{eqnarray*}
where $x=[a_{0},a_{1},a_{2},a_{3},...]$ stands for the representation of $x$ by
a (regular) continued fraction \cite{khin} (note that for $x\in[0,1)$ we have
$a_{0}=0$). This function was proved to be continuous, monotonic and singular
\cite{denjoy}. For convenience we introduce the Minkowski question mark function
$F(x)$ on the whole non-negative real line by
\begin{eqnarray*}
F([a_{0},a_{1},a_{2},a_{3},...])=1-2^{-a_{0}}+2^{-(a_{0}+a_{1})}-2^{-(a_{0}+a_{1}+a_{2})}+....
\label{min}
\end{eqnarray*}
Thus, $2F(x)=?(x)$ for $x\in[0,1]$. The function $F(x)$ as a distribution
function is uniquely determined by the functional equation \cite{ga}
\begin{eqnarray}
2F(x)=\left\{\begin{array}{c@{\qquad}l} F(x-1)+1 & \mbox{if}\quad x\geq 1,
\\ F({x\over 1-x}) & \mbox{if}\quad 0\leq x<1. \end{array}\right.
\label{distr}
\end{eqnarray}
This implies $F(x)+F(1/x)=1$.
\begin{figure}[h]
\begin{center}
\includegraphics[width=290pt,height=390pt,angle=-90]{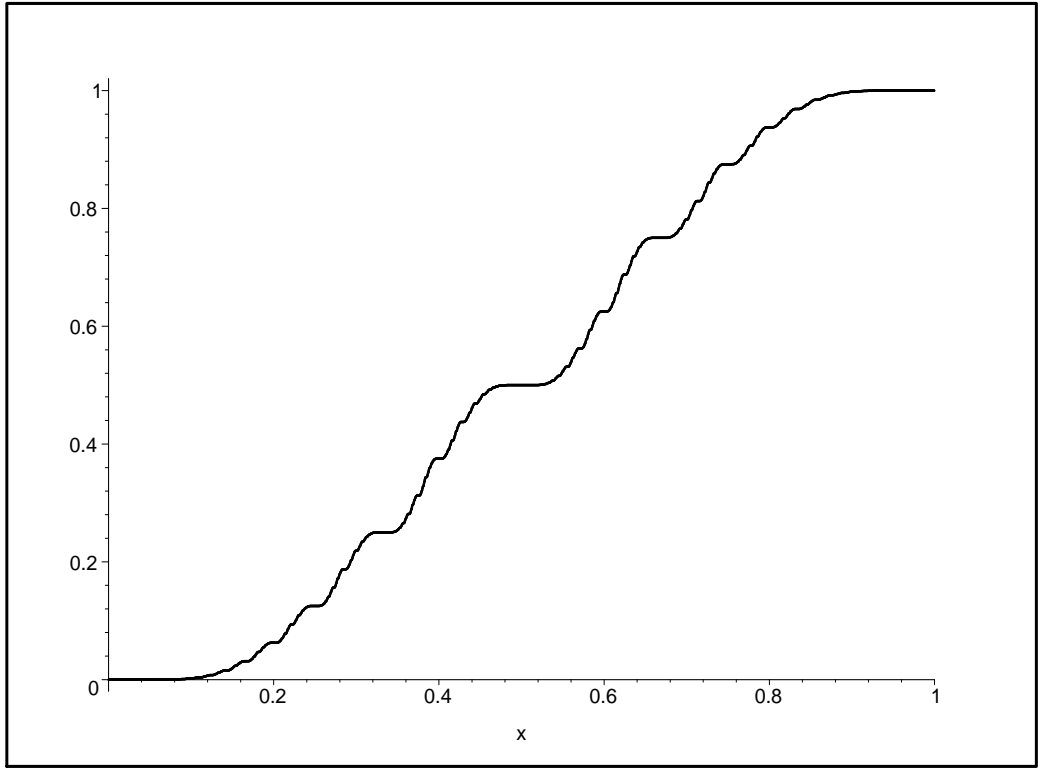}
\caption{Minkowski's question mark $?(x)$, $x\in[0,1]$}
\end{center}
\end{figure}
 The number-theoretic significance of the function
$F(x)$ is apparent from the fact that it is exactly the limit probabilistic
distribution of rationals in the $n$th generation of the Calkin-Wilf tree
\cite{cw} (which is just a permutation of the Stern-Brocot tree; the Farey tree is a
subtree of both). This tree is defined as follows. Starting from the root
$\frac{1}{1}$, each fraction generates two offsprings
\begin{eqnarray*}
{a\over b}\quad\mapsto\quad {a\over a+b}\ ,\quad {a+b\over b}.
\end{eqnarray*}
Elementary considerations show that this tree contains any positive rational
number once and only once in lowest terms \cite{cw}. First four iterations lead
to
\[
\xymatrix @R=.5pc @C=.5pc { & & & & & & & {1\over 1} & & & & & & & \\
& & & {1\over 2} \ar@{-}[urrrr] & & & & & & & & {2\over 1} \ar@{-}[ullll] & & & \\
& {1\over 3} \ar@{-}[urr] & & & & {3\over 2}\ar@{-}[ull] & & & & {2\over 3}\ar@{-}[urr] & & & & {3\over 1} \ar@{-}[ull] & \\
{1\over 4} \ar@{-}[ur] & & {4\over 3} \ar@{-}[ul] & & {3\over 5} \ar@{-}[ur] &
& {5\over 2} \ar@{-}[ul] & & {2\over 5} \ar@{-}[ur] & & {5\over 3} \ar@{-}[ul]
& & {3\over 4} \ar@{-}[ur] & & {4\over 1} \ar@{-}[ul] }
\]
It is important to note that the $n$th generation of the
Calkin-Wilf binary tree consists of exactly those rational numbers
whose elements of the continued fraction sum up to $n$. This fact
can be easily inherited directly from the definition. First, if
rational number $\frac{a}{b}$ is represented as a continued
fraction $[a_{0},a_{1},...,a_{r}]$, then the map
$\frac{a}{b}\rightarrow\frac{a+b}{b}$ maps $\frac{a}{b}$ to
$[a_{0}+1,a_{1}...,a_{r}]$. Second, the map
$\frac{a}{b}\rightarrow\frac{a}{a+b}$ maps $\frac{a}{b}$ to
$[0,a_{1}+1,...,a_{r}]$ in case $\frac{a}{b}<1$, and to
$[1,a_{0},a_{1},...,a_{r}]$ in case $\frac{a}{b}>1$. This is an
important fact which makes the investigations of rational numbers
according to their position in the Calkin-Wilf tree highly
motivated from the perspective of metric number theory and
dynamics of continued fractions. We refer the reader to papers
\cite{ga}, \cite{ga1} and \cite{ga2} for a wider overview of the properties of $?(x)$, and to a web page \cite{mink} for an exhaustive
bibliography list.\\

The work \cite{ga} by the author initiated study of the moments of the Minkowski
question mark function. Surprisingly, the only work where an analogous sequence of moments was considered before is
\cite{grabner_tichy}, where the authors consider the moments of a related singular function $F_{2}(x)$. The moments of the
question mark function are defined by
\begin{eqnarray}
m_{L}=\int\limits_{0}^{\infty}\Big{(}\frac{x}{x+1}\Big{)}^{L}\d
F(x)=2\int\limits_{0}^{1}x^{L}\d F(x)=\int\limits_{0}^{1}x^{L}\d
?(x).\label{moments}
\end{eqnarray}
This sequence is of definite number-theoretical significance because
\begin{eqnarray}
m_{L}=\lim_{n\rightarrow\infty}2^{2-n}\sum\limits_{a_{1}+a_{2}+...+a_{s}=n}[0,a_{1},a_{2},...,a_{s}]^{L},
\label{mom}
\end{eqnarray}
(the summation takes place over rational numbers presented as continued
fractions; thus, $s\geq 1$, $a_{i}\geq 1$ and $a_{s}\geq 2$). As an example of structural
properties of the sequence $m_{L}$ we mention the following result.
 Let
$\mathfrak{m}(t)=\sum\limits_{L=0}^{\infty}\frac{m_{L}t^{L}}{L!}$, which is an entire
function. One of the main results about $\mathfrak{m}(t)$ is that it is uniquely
determined by the regularity condition $\mathfrak{m}(-t)\ll e^{-\sqrt{t\log2}}$ as
$t\rightarrow\infty$, the boundary condition $\mathfrak{m}(0)=1$, and the integral
equation \cite{ga}
\begin{eqnarray*}
\mathfrak{m}(-s)=(2e^{s}-1)\int\limits_{0}^{\infty}\mathfrak{m}'(-t)J_{0}(2\sqrt{st})\d t,\quad
s\in\mathbb{R}_{+}.
\end{eqnarray*}
(Here $J_{0}(*)$ stands for the Bessel function
$J_{0}(z)=\frac{1}{\pi}\int_{0}^{\pi}\cos(z\sin x)\d x$). This reveals the
surprising fact that $\mathfrak{m}(t)$ can be considered as the dyadic analogue of entire
functions associated with Maass wave forms \cite{lewis}. Also, the generating
function $G(z)=\sum_{L=1}^{\infty}m_{L}z^{L-1}$ can be analytically continued
to the cut plane $\mathbb{C}\setminus (1,\infty)$ and it satisfies the three
term functional equation
\begin{eqnarray*}
\frac{1}{z}+\frac{1}{z^{2}}G\Big{(}\frac{1}{z}\Big{)}+2G(z+1)=G(z).
\end{eqnarray*}
Minding the specific form of this equation, it is natural to call $G(z)$
\it the dyadic period function. \rm \\

It is not clear whether there exists a closed form formula for the moments
$m_{L}$. This would be greatly desirable minding the expression (\ref{mom}). On
the other hand, one has the following asymptotic result \cite{ga1}. Let ${\sf
C}=e^{-2\sqrt{\log2}}=0.189169995269+$. Then the following estimate holds, as
$L\rightarrow\infty$:
\begin{eqnarray*}
{\sf C}^{\sqrt{L}}\ll m_{L}\ll L^{1/4}{\sf C}^{\sqrt{L}}.
\end{eqnarray*}
Both implied constants are absolute. The aim of this paper is to find the
first exact asymptotic term. Though from our point of view this is of inferior
significance than the closed form formula (if the latter does exist),
nevertheless, it is of definite interest too.

\begin{thm}Let $m_{L}$ be defined by (\ref{moments}). Then the
following holds as $L\rightarrow\infty$:
\begin{eqnarray*}
m_{L}=\sqrt[4]{4\pi^{2}\log2}\cdot c_{0}\cdot L^{1/4}{\sf C}^{\sqrt{L}}+
O(L^{-1/4}{\sf C}^{\sqrt{L}}).
\end{eqnarray*}
Here the constant $c_{0}$ is given by $c_{0}=\int_{0}^{1}\Psi(x)\d
x=1.030199563382+$, where $\Psi(x)$ is $1-$periodic function given by
$\Psi(x)=2^{x}(1-F(x))$ for $x\in[0,\infty)$.
\end{thm}

{\it Remark 1. } The periodicity of $\Psi(x)$ is obvious from (\ref{distr}).
Numerically, $\sqrt[4]{4\pi^{2}\log2}\cdot c_{0}=2.356229889908+$. The sequence of
moments
\begin{eqnarray*}
M_{L}=\int\limits_{0}^{\infty}x^{L}\d F(x)
\end{eqnarray*}
is equally important. These moments represent values of higher left
derivatives of $G(z)$ at $z=1$; more precisely, a formal Taylor expansion at
$z=1$ is given by
\begin{eqnarray*}
G(z+1)=\sum\limits_{L=1}^{\infty}M_{L}z^{L-1}.
\end{eqnarray*}
As it was proved in \cite{ga}, one has the asymptotic formula
$M_{L}\sim L!\frac{c_{0}}{(\log 2)^{L}}$. Moreover, in this case
there exists an exact convergent asymptotic series \cite{ga1}. Though
sequences $M_{L}$ and $m_{L}$ are linearly dependent via relations
\begin{eqnarray*}
m_{L}=M_{L}-\sum\limits_{s=0}^{L-1}M_{s}\binom{L}{s},\quad L\geq 0,
\end{eqnarray*}
it is nevertheless significant that the same structural constant $c_{0}$ does
manifest in asymptotic formulae for both of these sequences.\\
\indent {\it Remark 2. }As it is clear from the proof, other terms of
asymptotic expansion can be calculated as well. We confine to the
first term since calculations are standard (though tedious) with
no new ideas being introduced. Note that the structural constants
$c_{s}=\int_{0}^{1}x^{s}\Psi(x)\d x$ do appear in the asymptotic
expansion, and the method used shares some similarities with
the Euler-Maclaurin summation.\\
\indent {\it Remark 3. }If we start from the representation of $m_{L}$ via
the last integral of (\ref{moments}) (which is the most natural),
then to obtain the representation by the first integral we use the
second equation in (\ref{distr}). Moreover, one of the core stages
of the proof relies on the periodicity of $\Psi(x)$ and so uses
the first equation in (\ref{distr}). Since these two equations are
characteristic only to the Minkowski question mark function, generally
speaking, our asymptotic formula for the moments of $F(x)$ is
unique among similar results for other probabilistic distributions
with proper support on the interval $[0,1]$.
\section{Proof}
\subsection{Preliminary calculations}

In this paper we use the notation of Landau to denote by
$\mathcal{B}$ some absolutely bounded function in certain
neighborhood of a variable; in our case $\mathcal{B}$ depends
mostly on $L$ and we consider the case $L\rightarrow\infty$. As a
convention, $\mathcal{B}$ stands for different function if
considered in another or even the same formula. Occasionally, this
notation is used to denote dependence on other variables and in
each case it should be clear what variable and neighborhood is
implied in a notation $\mathcal{B}$. The main tool of the proof is a
Laplace method. The latter is used in finding
an asymptotic expansion of certain integrals depending on parameter. Though
in our case certain corrections, amendments and variations are
necessary, the main technique is standard and can be found, for
example, in \cite{lavr}. Since $1-F(x)=2^{-x}\Psi(x)$, identity
(\ref{moments}) implies
\begin{eqnarray*}
m_{L}&=&\int\limits_{0}^{\infty}\Big{(}\frac{x}{x+1}\Big{)}^{L}\d (F(x)-1)=
(F(x)-1)\Big{(}\frac{x}{x+1}\Big{)}^{L}\Big{|}_{0}^{\infty}\\&+&L\int\limits_{0}^{\infty}
2^{-x}\frac{x^{L-1}}{(x+1)^{L+1}}\Psi(x)\d x=L\int\limits_{0}^{\infty}
2^{-x}\frac{x^{L-1}}{(x+1)^{L+1}}\Psi(x)\d x.
\end{eqnarray*}
Put $c=\log 2$. The function $\mathbf{f}(x)=L\cdot\log\frac{x}{x+1}-cx$
achieves its maximum at $x=x_{0}$, where
\begin{eqnarray*}
x_{0}(x_{0}+1)=\frac{L}{c};\text{ thus, }
x_{0}=\frac{-1+\sqrt{1+\frac{4L}{c}}}{2},\text{ consequently }
x_{0}=\sqrt{\frac{L}{c}}-\frac{1}{2}+\frac{\mathcal{B}}{\sqrt{L}}.
\end{eqnarray*}
Let $c_{0}=\int_{0}^{1}\Psi(x)\d x$, and let us rewrite the expression for
$m_{L}$ as
\begin{eqnarray*}
m_{L}=Lc_{0}\int\limits_{0}^{\infty}\frac{1}{x(x+1)}e^{\mathbf{f}(x)}\d x+
L\int\limits_{0}^{\infty}(\Psi(x)-c_{0})\frac{1}{x(x+1)}e^{\mathbf{f}(x)}\d x=
c_{0}Lg_{L}+Lr_{L}.
\end{eqnarray*}
Our specific choice in extracting $c_{0}$ out of $\Psi(x)$ as a dominant
ingredient can be motivated for the following reason. Since
$0.9<\Psi(x)<1.2$, the main
weight of the integral defining $m_{L}$ (as far as $\Psi(x)$ is concerned)
befalls on a certain constant in the range $(0.9,1.2)$. Moreover, it is
easy to verify that for any continuous $1-$periodic function $\Upsilon(x)$ one has
\begin{eqnarray*}
\int\limits_{-\infty}^{\infty}\Upsilon(x)e^{-Ax^{2}}\d
x&=&\frac{\sqrt{\pi}\Upsilon(0)}{\sqrt{A}}+O(A^{-1})\text{ as
}A\rightarrow\infty\text{ (if }\Upsilon(x)\text{ is smooth)},\\
\int\limits_{-\infty}^{\infty}\Upsilon(x)e^{-Ax^{2}}\d
x&=&\frac{\sqrt{\pi}\int\limits_{0}^{1}\Upsilon(x)\d x}{\sqrt{A}}+O(1)\text{ as
}A\rightarrow 0+.
\end{eqnarray*}
This can be seen empirically from the fact that in the first case the dominant
weight of the integral is supported only in the neighborhood of $0$, while in
the second case this interval has a length tending to infinity. As it is seen
from the next subsection, we have the second case ({\it a posteriori}, the
weight of $m_{L}$ is properly supported on the interval
$(\sqrt{L/c}-L^{2/7},\sqrt{L/c}+L^{2/7})$, and hence the correct constant is
$c_{0}$ rather then $\Psi(x_{0})$.

\subsection{ Evaluation of $g_{L}$.}

 Choose $\delta=\delta(L)$ such that $\frac{\delta}{L^{1/4}}\rightarrow\infty$
and $\frac{\delta}{L^{1/3}}\rightarrow 0$. For example, henceforth we fix
$\delta=\delta(L)=L^{2/7}$. Separate the integral defining $g_{L}$ into four
parts:
\begin{eqnarray}
g_{L}=\int\limits_{0}^{1}+\int\limits_{1}^{x_{0}-\delta}+\int\limits_{x_{0}-\delta}^{x_{0}+\delta}+
\int\limits_{x_{0}+\delta}^{\infty}\frac{1}{x(x+1)}e^{\mathbf{f}(x)}\d
x=\mathscr{I}_{1}+\mathscr{I}_{2}+\mathscr{I}_{3}+\mathscr{I}_{4}.\label{int}
\end{eqnarray}
\subsubsection{Evaluation of $\mathscr{I}_{3}$}

 First, $\mathbf{f}'(x_{0})=0$ and $\mathbf{f}''(x_{0})<0$. Suppose
$|x-x_{0}|\leq\delta$. Then the Taylor formula gives that for certain
$\theta_{x}\in[x_{0},x]$ one has
\begin{eqnarray*}
\mathbf{f}(x)=\mathbf{f}(x_{0})-\alpha(x-x_{0})^{2}+\beta(x-x_{0})^{3}+\frac{f^{(4)}(\theta_{x})}{24}(x-x_{0})^{4}.
\end{eqnarray*}
Direct calculations show that
\begin{eqnarray*}
\mathbf{f}''(x)&=&-\frac{L(2x+1)}{x^{2}(x+1)^{2}};\text{ thus, }
\mathbf{f}''(x_{0})=-\frac{c^{2}}{L}(2x_{0}+1);\text{ this implies }
\alpha=\frac{c^{3/2}}{\sqrt{L}}+\frac{\mathcal{B}}{L^{3/2}};\\
\beta&=&\frac{\mathbf{f}'''(x_{0})}{6}=\frac{L(3x^{2}_{0}+3x_{0}+1)}{3x_{0}^{3}(x_{0}+1)^{3}}=
\frac{c^{2}}{L}+\frac{\mathcal{B}}{L^{3/2}};\\
\mathbf{f}^{(4)}(x)&=&6L\Big{(}\frac{1}{(x+1)^{4}}-\frac{1}{x^4}\Big{)};\text{ this gives }
\mathbf{f}^{(4)}(\theta_{x})=\frac{\mathcal{B}}{L^{3/2}}.
\end{eqnarray*}
In the same fashion,
\begin{eqnarray*}
\frac{1}{x(x+1)}=\frac{1}{x_{0}(x_{0}+1)}+\gamma(x-x_{0})+\sigma_{x}(x-x_{0})^{2},
\end{eqnarray*}
where
\begin{eqnarray*}
\gamma=-\frac{2x_{0}+1}{x^{2}_{0}(x_{0}+1)^{2}}=-\frac{2c^{3/2}}{L^{3/2}}+\frac{\mathcal{B}}{L^{5/2}},\quad
\sigma_{x}=\frac{\mathcal{B}}{L^{2}}.
\end{eqnarray*}
Now let us evaluate the value of $e^{\mathbf{f}(x_{0})}$. For
$x\rightarrow\infty$, we have
\begin{eqnarray*}
\log\frac{x}{x+1}&=&-\frac{1}{x}+\frac{1}{2x^{2}}+\frac{\mathcal{B}}{x^{3}};\text{ this gives}\\
L\log\frac{x_{0}}{x_{0}+1}&=&-\frac{L}{\sqrt{\frac{L}{c}}-\frac{1}{2}+\frac{\mathcal{B}}{\sqrt{L}}}+\frac{c}{2}+\frac{\mathcal{B}}{\sqrt{L}}=
-\sqrt{cL}+\frac{\mathcal{B}}{\sqrt{L}}.
\end{eqnarray*}
Further,
\begin{eqnarray*}
-cx_{0}&=&-\sqrt{cL}+\frac{c}{2}+\frac{\mathcal{B}}{\sqrt{L}};\text{ this gives}\\
e^{\mathbf{f}(x_{0})}&=&\exp\Big{(}-2\sqrt{cL}+\frac{c}{2}+\frac{\mathcal{B}}{\sqrt{L}}\Big{)}=
e^{-2\sqrt{cL}}\sqrt{2}\Big{(}1+\frac{\mathcal{B}}{\sqrt{L}}\Big{)}.
\end{eqnarray*}
We already obtained all necessary components to evaluate $\mathscr{I}_{3}$.
Note that
\begin{eqnarray*}
\beta(x-x_{0})^{3}=\frac{\mathcal{B}}{L}\cdot L^{6/7}=o(1),\quad
\frac{\mathcal{B}}{L^{3/2}}(x-x_{0})^{4}=o(1).
\end{eqnarray*}
Since $e^{y}=1+y+\mathcal{B}y^{2}$ for $y=o(1)$, the function under the integral (\ref{int}) in
the range $(x_{0}-\delta,x_{0}+\delta)$ can be written as
\begin{eqnarray*}
\frac{1}{x(x+1)}e^{\mathbf{f}(x)}&=&\Big{(}\frac{1}{x_{0}(x_{0}+1)}+\gamma(x-x_{0})+\frac{\mathcal{B}}{L^{2}}(x-x_{0})^{2}\Big{)}\\
&\times&\exp\Big{(}\mathbf{f}(x_{0})-\alpha(x-x_{0})^{2}+\beta(x-x_{0})^{3}+\frac{\mathcal{B}}{L^{3/2}}(x-x_{0})^{4}\Big{)}\\
&=&e^{\mathbf{f}(x_{0})}e^{-\alpha(x-x_{0})^{2}}\cdot\Big{(}1+\beta(x-x_{0})^{3}+\frac{\mathcal{B}}{L^{3/2}}(x-x_{0})^{4}
+\frac{\mathcal{B}}{L^{2}}(x-x_{0})^{6}\Big{)}\\
&\times&\Big{(}\frac{1}{x_{0}(x_{0}+1)}+\gamma(x-x_{0})+\frac{\mathcal{B}}{L^{2}}(x-x_{0})^{2}\Big{)}
\end{eqnarray*}
(The bound $|x-x_{0}|=\mathcal{B}L^{1/2}$ was used in merging
$\frac{\mathcal{B}\beta}{L^{3/2}}(x-x_{0})^{7}$ and
$\frac{\mathcal{B}}{L^{3}}(x-x_{0})^{8}$ into
$\frac{\mathcal{B}}{L^{2}}(x-x_{0})^{6}$). Thus, we have decomposed the
function under integral (\ref{int}) into the sum of twelve functions. It is important to
note that two of these functions
\begin{eqnarray*}
e^{\mathbf{f}(x_{0})}e^{-\alpha(x-x_{0})^{2}}\gamma(x-x_{0})\text{ and }
e^{\mathbf{f}(x_{0})}e^{-\alpha(x-x_{0})^{2}}\frac{1}{x_{0}(x_{0}+1)}\beta(x-x_{0})^{3},
\end{eqnarray*}
though being comparatively large, are odd function in $(x-x_{0})$ and
consequently contribute $0$ after integration. We are left with evaluating all
the rest summands. This is a routine job. For example,
\begin{eqnarray*}
&&e^{\mathbf{f}(x_{0})}\gamma\beta\int\limits_{x_{0}-\delta}^{x_{0}+\delta}e^{-\alpha(x-x_{0})^{2}}(x-x_{0})^{4}\d
x\\
&=&e^{\mathbf{f}(x_{0})}\gamma\beta\frac{1}{\alpha^{5/2}}\int\limits_{-\sqrt{\alpha}\delta}^{\sqrt{\alpha}\delta}y^{4}e^{-y^{2}}\d
y=\mathcal{B}\frac{L^{5/4}}{L^{5/2}}{\sf
C}^{\sqrt{L}}=\mathcal{B}L^{-5/4}{\sf C}^{\sqrt{L}}.
\end{eqnarray*}
The second example:
\begin{eqnarray*}
& &e^{\mathbf{f}(x_{0})}\frac{\mathcal{B}\gamma}{L^{3/2}}\int\limits_{x_{0}-\delta}^{x_{0}+\delta}e^{-\alpha(x-x_{0})^{2}}|x-x_{0}|^{5}\d
x\\
&=&e^{\mathbf{f}(x_{0})}\frac{\mathcal{B}\gamma}{L^{3/2}\alpha^{3}}\int\limits_{-\sqrt{\alpha}\delta}^{\sqrt{\alpha}\delta}|y|^{5}e^{-y^{2}}\d
y=\mathcal{B}L^{-3/2-3/2+3/2}{\sf C}^{\sqrt{L}}=\mathcal{B}L^{-3/2}{\sf
C}^{\sqrt{L}}.
\end{eqnarray*}
Another example:
\begin{eqnarray*}
e^{\mathbf{f}(x_{0})}\frac{\mathcal{B}}{x_{0}(x_{0}+1)L^{3/2}}\int\limits_{x_{0}-\delta}^{x_{0}+\delta}e^{-\alpha(x-x_{0})^{2}}(x-x_{0})^{4}\d
x=\mathcal{B}\frac{L^{5/4}}{L^{5/2}}{\sf
C}^{\sqrt{L}}&=&\mathcal{B}L^{-5/4}{\sf C}^{\sqrt{L}}.
\end{eqnarray*}
 Therefore, a direct inspection shows that all functions apart
from the main term contribute at most $\mathcal{B}{\sf C}^{\sqrt{L}}L^{-5/4}$
into the value of $g_{L}$. We are left with determining the magnitude of the
main term. Since $\sqrt{\alpha}\delta\sim c^{3/4}L^{1/28}$, this yields
\begin{eqnarray*}
\int\limits_{-\sqrt{\alpha}\delta}^{\sqrt{\alpha}\delta}e^{-y^{2}}\d y=
\int\limits_{-\infty}^{\infty}e^{-y^{2}}\d
y+\mathcal{B}\int\limits_{\sqrt{\alpha}\delta}^{\infty}e^{-y^{2}}\d
y=\sqrt{\pi}+\mathcal{B}\exp(-c^{3/2}L^{1/14}).
\end{eqnarray*}
And so, the main term is given by
\begin{eqnarray*}
\frac{1}{x_{0}(x_{0}+1)}e^{\mathbf{f}(x_{0})}\int\limits_{-\delta}^{\delta}e^{-\alpha
x^{2}}\d
x&=&\Big{(}\sqrt{\pi}+\mathcal{B}\exp(-c^{3/2}L^{1/14})\Big{)}\frac{1}{\sqrt{\alpha}}\frac{1}{x_{0}(x_{0}+1)}e^{\mathbf{f}(x_{0})}\\
&=&\frac{c}{L}\sqrt{2\pi}e^{-2\sqrt{cL}}
\sqrt[4]{L}c^{-3/4}\Big{(}1+\frac{\mathcal{B}}{\sqrt{L}}\Big{)}.
\end{eqnarray*}
Thus, we obtain the main asymptotic term
\begin{eqnarray*}
\mathscr{I}_{3}=\sqrt[4]{4\pi^{2}\log 2}\cdot L^{-3/4}{\sf
C}^{\sqrt{L}}+\mathcal{B}L^{-5/4}{\sf C}^{\sqrt{L}}.
\end{eqnarray*}
\subsubsection{Evaluation of $\mathscr{I}_{1}$, $\mathscr{I}_{2}$ and $\mathscr{I}_{4}$}

Trivially,
\begin{eqnarray*}
\mathscr{I}_{1}=\int\limits_{0}^{1}\frac{x^{L-1}}{(x+1)^{L+1}}2^{-x}\d
x=\mathcal{B}2^{-L}.
\end{eqnarray*}

To evaluate $\mathscr{I}_{2}$, we use a standard inequality
\begin{eqnarray}
\log\Big{(}\frac{x}{x+1}\Big{)}<-\frac{1}{x}+\frac{1}{2x^{2}}\text{ for }x\geq
1.\label{nell}
\end{eqnarray}
Thus, since $\mathbf{f}(x)$ is an increasing function in the interval $[1,x_{0}]$,
we have
\begin{eqnarray*}
\mathscr{I}_{2}=\int\limits_{1}^{x_{0}-\delta}
\frac{1}{x(x+1)}e^{\mathbf{f}(x)}\d x<\log 2\cdot e^{\mathbf{f}(x_{0}-\delta)}=
\mathcal{B}e^{\mathbf{f}(x_{0}-\delta)}.
\end{eqnarray*}
Further, using the inequality (\ref{nell}), we have
\begin{eqnarray*}
\mathbf{f}(x_{0}-\delta)&=&L\log\Big{(}\frac{x_{0}-\delta}{x_{0}-\delta+1}\Big{)}-c(x_{0}-\delta)<
-\frac{L}{x_{0}-\delta}+\frac{L}{2(x_{0}-\delta)^{2}}-cx_{0}+c\delta\\
&=&-\frac{L}{x_{0}}-cx_{0}+\mathcal{B}+c\delta+\Big{(}\frac{L}{x_{0}}-\frac{L}{x_{0}-\delta}\Big{)}\\
&=&-\frac{L}{x_{0}}-cx_{0}+\mathcal{B}+c\delta-\frac{L\delta}{x_{0}(x_{0}+1)}-L\delta
\Big{(}\frac{1}{x_{0}(x_{0}-\delta)}-\frac{1}{x_{0}(x_{0}+1)}\Big{)}\\
&=&-\frac{L}{x_{0}}-cx_{0}+\mathcal{B}-\frac{L\delta^{2}}{x_{0}^{3}}=-2\sqrt{cL}+\mathcal{B}-c^{3/2}L^{1/14}
\end{eqnarray*}
(recall that $\delta=L^{2/7}$). Therefore,
\begin{eqnarray*}
\mathscr{I}_{2}=\mathcal{B}{\sf
C}^{\sqrt{L}}\exp(-c^{3/2}L^{1/14})=\mathcal{B}L^{-5/4}{\sf C}^{\sqrt{L}}.
\end{eqnarray*}
In the same vein,
\begin{eqnarray*}
\mathbf{f}(x_{0}+\delta)=-2\sqrt{cL}+\mathcal{B}-c^{3/2}L^{1/14}.
\end{eqnarray*}
Since $\mathbf{f}(x)$ is decreasing in the interval $[x_{0},\infty)$, we obtain
\begin{eqnarray*}
\mathscr{I}_{4}&=& \int\limits_{x_{0}+\delta}^{\infty}\frac{1}{x(x+1)}e^{\mathbf{f}(x)}\d
x<e^{\mathbf{f}(x_{0}+\delta)}\int\limits_{1}^{\infty}\frac{\d x}{x(x+1)}\\
&=&\mathcal{B}{\sf C}^{\sqrt{L}}\exp(-c^{3/2}L^{1/14})=\mathcal{B}L^{-5/4}{\sf
C}^{\sqrt{L}}.
\end{eqnarray*}
Combining all the estimates for $\mathscr{I}_{1}$, $\mathscr{I}_{2}$,
$\mathscr{I}_{3}$ and $\mathscr{I}_{4}$, we eventually obtain
\begin{eqnarray*}
g_{L}=\sqrt[4]{4\pi^{2}\log 2}\cdot L^{-3/4}{\sf C}^{\sqrt{L}}+\mathcal{B}L^{-5/4}{\sf
C}^{\sqrt{L}}.
\end{eqnarray*}
\subsection{Evaluation of $r_{L}$.}

In this subsection we can be more concise, since the method is the same as we
used to evaluate $g_{L}$. Recall that
\begin{eqnarray*}
r_{L}=\int\limits_{0}^{\infty}(\Psi(x)-c_{0})\frac{1}{x(x+1)}e^{\mathbf{f}(x)}\d
x\end{eqnarray*} Let
\begin{eqnarray*}
\hat{\Psi}(x)=\int\limits_{0}^{x}(\Psi(t)-c_{0})\d t.
\end{eqnarray*}
The definition of $c_{0}$ yields that $\hat{\Psi}(x)$ is bounded (and periodic)
function. Using an integration by parts, we can rewrite an expression for $r_{L}$ as
\begin{eqnarray*}
r_{L}=\int\limits_{0}^{\infty}\hat{\Psi}(x)h(x)e^{\mathbf{f}(x)}\d x; \text{ thus, }
r_{L}=\mathcal{B}\int\limits_{0}^{\infty}h(x)e^{\mathbf{f}(x)}\d x,
\end{eqnarray*}
where
\begin{eqnarray*}
h(x)=\frac{1}{x^{2}}-\frac{1}{(x+1)^{2}}-\frac{1}{x(x+1)}\mathbf{f}'(x)=
-\frac{L}{x^{2}(x+1)^{2}}+\frac{c}{x(x+1)}+\frac{2x+1}{x^{2}(x+1)^{2}}.
\end{eqnarray*}
Note that
\begin{eqnarray*}
h(x_{0})=\frac{2x_{0}+1}{x^{2}_{0}(x_{0}+1)^{2}}=\frac{2c^{3/2}}{L^{3/2}}+\frac{\mathcal{B}}{L^{5/2}}.
\end{eqnarray*}
As a matter of fact, we have a completely analogous integral to the one defining
$g_{L}$, with the major difference in the estimate
$h(x_{0})=\mathcal{B}L^{-3/2}$, whereas
$\frac{1}{x_{0}(x_{0}+1)}=\mathcal{B}L^{-1}$. Thus, using the same method to
evaluate $r_{L}$ as we did with $g_{L}$, one gets
\begin{eqnarray*}
r_{L}=\mathcal{B}L^{1/4-3/2}{\sf C}^{\sqrt{L}}=\mathcal{B}L^{-5/4}{\sf
C}^{\sqrt{L}}.
\end{eqnarray*}
Since $m_{L}=c_{0}Lg_{L}+Lr_{L}$, this finishes the proof of Theorem 1.
$\blacksquare$
\subsection{Final remarks}
As can be inherited from the proof, the share of integrals $\mathscr{I}_{1}$,
$\mathscr{I}_{2}$ and $\mathscr{I}_{4}$ into the exact value of $m_{L}$ is of
order ${\sf C}^{\sqrt{L}}\exp(-c^{3/2}L^{1/14})$. Thus, a standard machinery of
an asymptotic expansion shows that all terms for asymptotic series of $m_{L}$ are
supported on the neighborhood $(x_{0}-\delta,x_{0}+\delta)$. More thorough
inspection reveals that
\begin{eqnarray*}
m_{L}\sim L^{1/4}{\sf
C}^{\sqrt{L}}\cdot\sum\limits_{i=0}^{\infty}\frac{{\sf
A}_{i}}{L^{i/2}}.
\end{eqnarray*}
As was mentioned, every term of this expansion can be obtained by
(increasingly tedious) computations. For example, ${\sf A}_{1}$ is
a sum of $c_{0}$ and $c_{1}=\int_{0}^{1}x\Psi(x)\d x$, each
multiplied by some explicit constant (constant $c_{1}$ occurs
while extracting the main term of the integral $r_{L}$). Nevertheless, from a number-theoretic
point of view our main concern is the structure of values $m_{L}$ rather than
their magnitude (though dominant) in some small neighborhood of $x_{0}$, and
hence we have confined only in
exhibiting such possibility of asymptotic expansion.\\

We finish with providing a table for some values of constants $m_{L}$.

\begin{center}
\begin{tabular}{|r|l|l|}
\hline
\multicolumn{3}{|c|}{\textbf{Sequence $m_{L}$}}\\
\hline
$L$ & $m_{L}$& $m^{\star}_{L}$\\
\hline
$1$ & $0.5000000000$    &$2.643125297$\\
$2$ & $0.2909264764$    &$2.577573745$\\
$3$ & $0.1863897146$    &$2.533204605$\\
$4$ & $0.1269922584$    &$2.509329792$\\
$5$ & $0.09016445494$    &$2.496320715$\\
$6$ & $0.06592816257$    &$2.488147649$\\
$7$ & $0.04929431046$    &$2.481940613$\\
$8$ & $0.03751871185$   &$2.476544438$\\
$9$ & $0.02897962203$   &$2.471583746$\\
$10$ & $0.02266585817$  &$2.466982861$\\
$11$ & $0.01792085923$  &$2.462750421$\\
$12$ & $0.01430468951$  &$2.458897371$\\
$20$ & $0.003008686707$ &$2.438565967$\\
$30$ & $0.0006211064464$ &$2.425096683$\\
$40$ & $0.0001622371309$ &$2.416702495$\\
$50$ & $0.00004937221843$ &$2.410831724$\\
$100$ & $0.0000004445933003$ &$2.395743861$\\
\hline
\end{tabular}
\end{center}

\indent Here $m^{\star}_{L}=\frac{m_{L}}{\sqrt[4]{L}{\sf C}^{\sqrt{L}}}$.
The main result of the paper implies that the sequence $m^{\star}_{L}$
tends to the limit $c_{0}\sqrt[4]{4\pi^{2}\log 2}=2.3562298899+$.
Unfortunately, we do not have yet any evidence that the closed
form expression exists for $c_{0}$ (as well as for $m_{L}$ with
$L\geq 2$). Finally, we remark that the result of Theorem 1 should
be considered in conjunction with linear relations \cite{ga} $
m_{L}=\sum_{s=0}^{L}\binom{L}{s}(-1)^{s}m_{s}$, $L\geq 0$, which
the sequence $m_{L}$ satisfies.

\par\bigskip

\noindent

\noindent The School of Mathematical Sciences, The University of Nottingham, University Park, Nottingham NG7 2RD United Kingdom\\

\noindent Vilnius University, The Department of Mathematics and Informatics, Naugarduko 24, Vilnius, Lithuania\\
{\tt giedrius.alkauskas@maths.nottingham.ac.uk}\\

\smallskip

\begin{thebibliography}{9}

\bibitem{ga} {\sc G. Alkauskas}, The moments of Minkowski question mark function: the dyadic period function
(submitted); {\tt arXiv:0801.0051}.

\bibitem{ga1} {\sc G. Alkauskas}, Generating and zeta functions, structure,
spectral and analytic properties of the moments of Minkowski
question mark function (submitted); {\tt
arXiv:0801.0056}.

\bibitem{ga2}{\sc G. Alkauskas}, Minkowski question mark function and its generalizations, associated with {\bf p}-continued fractions:
fractals, explicit series for the dyadic period function and moments (submitted); {\tt arXiv:0805.1717}.

\bibitem{cw}{\sc N. Calkin, H. Wilf}, Recounting the rationals,
{\it Amer. Math. Monthly}, {\bf 107}:360-363, 2000.

\bibitem{denjoy} {\sc A. Denjoy}, Sur une fonction r\'{e}elle de Minkowski, {\it J. Math. Pures Appl.}, {\bf 17}:105-151, 1938.

\bibitem{grabner_tichy}{\sc P.J. Grabner, P. Kirschenhofer, R.F. Tichy},
    Combinatorial and arithmetical properties of linear numeration systems,
      {\it Combinatorica}, {\bf 22}(2):245-267, 2002.

\bibitem{khin}{\sc A.Ya. Khinchin}, {\it Continued fractions}, The University of Chicago Press, 1964

\bibitem{lavr}{\sc M.A. Lavrent'ev, B.V. Shabat},
{\it Methods of the theory of functions of a complex variable},
Nauka, Moscow, 1987 (in Russian).

\bibitem{lewis}{\sc J.B. Lewis}, Spaces of holomorphic functions equivalent to the even Maass cusp forms,
{\it Invent. Math.}, {\bf 127}(2):271-306, 1997.

\bibitem{mink} An exhaustive bibliography on the Minkowski question mark function,
\url{http://www.maths.nottingham.ac.uk/personal/pmxga2/minkowski.htm}
\end{thebibliography}
\end{document}